\documentclass{amsart}
\usepackage{amsmath, amsthm, amscd, amsfonts,graphicx,needspace,url,amssymb}
\theoremstyle{definition}
\newtheorem{Q}{Q\!}
\newtheorem{KQ}[Q]{(Q\!}

\renewcommand{\r}{\mathrm}

\begin{document}

\begin{center}
\texttt{Comments, corrections,
and related references welcomed, as always!}\\[.5em]
{\TeX}ed \today
\vspace{2em}
\end{center}

\title[Possible questions for the next Kourovka notebook]%
{Some questions for possible submission to the next Kourovka notebook}
\thanks{This note can be viewed at
\url{http://arxiv.org/abs/1904.04298} and at
\url{http://math.berkeley.edu/~gbergman/papers/unpub/418.preKourovka.pdf}\,.
The latter file may be updated more often than the former.}

\author{George M. Bergman}
\address{George M. Bergman, Department of Mathematics,
University of California, Berkeley, CA 94720-3840, USA}
\email{gbergman@math.berkeley.edu}
\urladdr{https://math.berkeley.edu/~gbergman}

\maketitle

\noindent{\bf Background.}

Some years ago, I had meant to submit a few questions to the 19th
edition of the Kourovka Notebook \cite{Kourovka}, but I was busy,
and it came out before I had done so.
So I decided to gather ahead of time questions I would
like to submit to the next edition.
Looking over past papers, I found
group-theoretic questions that I had raised, and
listed most of them, together with a few new ones.

I submitted the list; and several of the questions were
included in the 20th edition of the Kourovka Notebook,
and several more in the 21st.
Below, I have replaced those that were used
in those two editions with references to the questions there.
(I could have simply deleted them, but that would have meant
that the ones that followed would change their numbering,
which could complicate things for someone who had taken notes
regarding a past version of this list, and then looked at the
current version.)

I may add more questions in the future.

If you look this over,
let me know if you have any comments -- in particular if, to
your knowledge, any of the questions below have been asked before,
or if you know or can see any answers.

\vspace{1em}
\noindent{\bf Questions.}

\begin{KQ}\label{Q.n-u.p.} 
Appears as Question 20.12 in the Kourovka Notebook.)
%
%
\end{KQ}

\begin{Q}\label{Q.SLIO} 
If $S$ is a subset of a group $G,$ then an
element $s\in S$ is called an {\em extremal point} of $S$
if for every $g\in G-\{1\}$, either $sg^{-1}\notin S$ or $sg\notin S$.
A group $G$ is said to be {\em diffuse} if
every finite nonempty subset of $G$ has an extremal point.
In \cite[Proposition~6.2]{PL+DWM}, this condition on $G$
is shown equivalent (inter alia) to the condition called LIO (locally
invariant order), namely that $G$ admits a set-theoretic
total ordering $\geq$ such that for all $s\in G$
and \,$g\in G-\{0\},$ at least one of $sg^{-1}>s$ or $sg>s$ holds.

Let us say that $G$ has SLIO (strongly locally invariant order)
if it admits an ordering such that for all such $s$ and $g$,
{\em exactly} one of $sg^{-1}>s$ or $sg>s$ holds.
(In other words, every sequence of the form $(sg^i)_{i\in\mathbb{Z}}$
is monotone increasing or monotone decreasing.)

It is not hard to show that
\needspace{2em}

\hspace{2em}{right-orderable} $\implies$ SLIO $\implies$ LIO
(equivalently, diffuse).\\
The composite implication is known not
to be reversible \cite[Appendix~B]{SK+JR}.

But is one or the other of the above two implications reversible?
\end{Q}

\begin{Q}\label{Q.nonunit} 
(i)~~\cite{overflow},
\cite[Definition~1 and Question~2]{prin_ids_in_kG} \ %
Let us call group $G$ {\em resistant} if for every
field $k,$ and every element $r=\sum_i c_i\,g_i$ of the group
algebra $k\,G$ whose support, ${\rm supp}(r)=\{g_i\mid c_i\neq 0\},$
has cardinality $>1,$ the $2$-sided ideal of $k\,G$
generated by $r$ is proper.

(Part (i) of this question, ``Are all free groups resistant?''
appears, reworded, as Question~20.9 in the Kourovka Notebook.
The remaining parts of this question, shown below, were not included
there.)
\vspace{.3em}

The question of whether all free groups are resistant
is the motivation of \cite{prin_ids_in_kG},
where many sorts of groups are proved non-resistant; e.g., any
group containing a nonabelian solvable subgroup.
Six more questions are asked there, of which I record two here.
\vspace{.3em}

(ii)~\cite[Question~19]{prin_ids_in_kG} \ %
Is the class of resistant groups closed under taking direct
products?\vspace{.3em}

(iii)~\cite[Definition~22 and Question~23]{prin_ids_in_kG} \ %
(Possible Freiheitssatz for group algebras) \ %
For $F$ the free group on generators $(x_i)_{i\in I},$ and $k$ a field,
let us call an element $r\in k\,F$ {\em strongly reduced} if, when
the elements of ${\rm supp}(r)$ are written as reduced words
in the free generators,
(a)~there is no symbol $x_i^{\pm 1}$ with which all
these elements begin,
(b)~there is no symbol $x_i^{\pm 1}$ with which all
these elements end, and
(c)~if $1\in\r{supp}(r)$ (so that the preceding
two conditions hold trivially)
there is no symbol $x_i^{\pm 1}$ such that all elements
of $\r{supp}(r)$ other than $1$ both begin with
$x_i^{\pm 1}$ and end with the inverse symbol, $x_i^{\mp 1}.$

For any $r\in k\,F,$ we shall say that a generator $x_{i_0}$ is
``involved in'' $r$ if $x_{i_0}$ or $x_{i_0}^{-1}$ occurs
anywhere in the normal form of any of the elements
of $\r{supp}(r).$

Suppose $I$ is the $2$-sided ideal of $k\,F$ generated by
a strongly reduced element $r.$
If $F'\subset F$ is the subgroup generated by $(x_i)_{i\in I-i_0},$
where $x_{i_0}$ is a generator involved in $r,$
must the composite map
$k\,F'\hookrightarrow k\,F\to (k\,F)/I$ be an embedding?
(A positive answer would imply a positive answer to~(i).)
\end{Q}

\begin{Q}\label{Q.fr_kG_in_D} 
\cite[Question 13]{kG_in_D}
Let $(G,\leq)$ be a right-ordered group, and $k$ a field.
Thus, the set of formal infinite sums $\sum_{g\in G} a_g g$
with coefficients $a_g\in k$ whose support
$\{g\in G\ |\ a_g\neq 0\}$ is well-ordered forms a right module
$k((G))$ over the group algebra $kG$, though this will not in
general have a natural structure of ring, or of left $kG$-module.

Dubrovin~\cite{Dubrovin} showed that every nonzero element
of $kG$ acts {\em invertibly} on the module $k((G))$.
We ask:

For each $x_1,\,x_2,\,y_1,\,y_2\in kG-\{0\},$ will the
right action of
$y_1 y_2^{-1} - x_1 x_2^{-1}$ on $k((G))$ either be invertible or zero?
(Here $y_2^{-1}$ and $x_2^{-1}$ denote the inverses of
the {\em actions} of $y_2$ and $x_2$.)

The above is easily seen to be equivalent to the same question
for the actions of $x_1^{-1} y_1 - x_2^{-1} y_2,$ of
$y_1 - x_1 x_2^{-1} y_2,$ and of $1 - x_1 x_2^{-1} y_2\,y_1^{-1}$.

An affirmative answer to this question
is equivalent to a matrix-theoretic property
asked for in \cite[Question~12]{kG_in_D}.
Such a result would be the ``next step'', after Dubrovin's result,
toward proving that $kG$ is embeddable in a division ring,
and indeed generates a division ring of operators on $k((G))$.
A countable chain of conditions, having Dubrovin's result as its
first step and the property asked for above as its second,
and which all together are equivalent to generating a division ring
of operators, is developed in \cite[\S10]{kG_in_D}.

For simplicity I am posing this question for $k$
a field rather than a general division ring.
If a positive answer is obtained, one should, of course, see
whether commutativity of $k$ is really needed, and whether
the result can be extended to group rings twisted by actions of
$G$ on $k$.
\end{Q}


\begin{KQ}\label{Q.B,B'} 
Appears as Question 20.10 in the Kourovka Notebook.)
%
%
\end{KQ}

\begin{KQ}\label{Q.B_gp+} 
Appears as Question 21.11 in the Kourovka Notebook.)
%
\end{KQ}

\begin{KQ}\label{Q.ultra} 
Appears as Question 21.12 in the Kourovka Notebook.)
%
\end{KQ}

\begin{KQ}\label{Q.ab<gp} 
Appears as Question 20.13 in the Kourovka Notebook.)
\end{KQ}

\begin{KQ}\label{Q.free_dual} 
Appears as Question 21.13 in the Kourovka Notebook.)
\end{KQ}

\begin{KQ}\label{Q.fr_inner} 
Appears as Question 21.14 in the Kourovka Notebook.)
%
\end{KQ}

\begin{Q}\label{Q.XYZ} 
\cite[Question 16, p.157]{nomial} \ %
Suppose $G$ is a group, and $X,\ Y,\ Z$ are finite
transitive $G$-sets, such that

\hspace{2em}$\r{g.c.d.}(\r{card}(X),\ \r{card}(Y),
\ \r{card}(Z))\ =\ 1,$\\
but such that none of the $G$-sets
$X\!\times\!Y,$ \ $Y\!\times\!Z,$ \ $Z\!\times\!X$
is transitive (so that
no two of $\r{card}(X)$,\, $\r{card}(Y)$,\, $\r{card}(Z)$
are relatively prime).

What, if anything, can one conclude about the orbit-structure
of the $G$-set $X\times Y\times Z$?
(For the motivation, see
\cite{nomial}, paragraph preceding this question, then
the group-theoretic proof of \cite[Theorem~1, p.139]{nomial},
and finally, the first paragraph of \S2 on that page.
But perhaps the question is, nonetheless, too vague
to include in the {\em Kourovka Notebook}.)
\end{Q}

\begin{KQ}\label{Q.cP_w_amalg} 
Appears as Question 21.15 in the Kourovka Notebook.)
%
\end{KQ}

\begin{KQ}\label{Q.2_vs_3} 
Appears as Question 20.14 in the Kourovka Notebook.)
%
%
%
\end{KQ}

\begin{KQ}\label{Q.resid_fin} 
Appears as Question 20.15 in the Kourovka Notebook.)
%
\end{KQ}

\begin{KQ}\label{Q.width_as_monoid} 
Appears as Question 21.16 in the Kourovka Notebook.)
%
\end{KQ}

\begin{Q}\label{Q.num_gens,rels} 
\cite[Question 10, p.437]{inf_fr_fin} \ %
In \cite{B+S.af} Baumslag and Shalen show that
a group is infinite if it has a presentation in which all
relations have the form $s^n = 1$,
for various group words $s$, and a common exponent $n$ which
is a prime power, such that the number of generators in the presentation
is at least 1 more than $1/n$ times the number of relations.

Can one get a more general criterion for a group to be
infinite, again based on counting
a relation of the form $s^n = 1$ as ``$1/n$ of a relation''
but with the ``common exponent'' assumption
replaced by the condition that the exponents $n$
so treated all be powers of a common prime $p$; or alternatively,
with the common exponent requirement retained,
but no requirement that it be a prime power;
or, covering both these generalizations, merely the
requirement that the set
of exponents so treated be totally ordered under divisibility?

(For the need for some such restrictions, see
\cite[Introduction]{inf_fr_fin}.
For another known result of the same sort, see \cite{Thomas},
also summarized in \cite[Introduction]{inf_fr_fin}.)
\end{Q}

\begin{KQ}\label{Q.pseudo_epis} 
Appears as Question 20.8 in the Kourovka Notebook.)
%
\end{KQ}

\begin{Q}\label{Q.HSP_f} 
(Part (i) below appears as Question 21.17 in the Kourovka Notebook.)

If $\mathbf{X}$ is a class of groups,
let $\mathbf{H}(\mathbf{X})$ denote the class of homomorphic
images of groups in $\mathbf{X}$,
let $\mathbf{S}(\mathbf{X})$ denote the class of groups isomorphic
to subgroups of groups in $\mathbf{X}$,
let $\mathbf{P}(\mathbf{X})$ denote the class of groups isomorphic
to products of families of groups in $\mathbf{X}$, and
let $\mathbf{P}_{\!f}(\mathbf{X})$ denote the class of groups isomorphic
to products of {\em finite} families of groups in $\mathbf{X}$.

By the group case of Birkhoff's Theorem,
$\mathbf{H\,S\,P}(\mathbf{X})$ is the variety
of groups generated by $\mathbf{X}$.

(i)  \cite[Question~27, p.281]{SHPS} \ %
If $\mathbf{X}$ is a class of {\em metabelian} groups, must
$\mathbf{H\,S\,P}_{\!f}(\mathbf{X})\subseteq
\mathbf{S\,H\,P\,S}(\mathbf{X})$?
(See preceding paragraph of \cite{SHPS} for motivation.)

The operators $\mathbf{H}$, $\mathbf{S}$, $\mathbf{P}$
make sense for more general algebraic structures,
in particular, lattice-ordered groups, as in

(ii)  (S.\,Comer, personal correspondence, cited in
\cite[Question~28(ii), p.281]{SHPS}.)
If $\mathbf{X}$ is a class of {\em lattice-ordered groups}, must
$\mathbf{S\,H\,P\,S}(\mathbf{X})=\mathbf{H\,S\,P}(\mathbf{X})$\,?
\end{Q}

\begin{Q}\label{Q.hyperid} 
A {\em hyperidentity} for a class $\mathbf{X}$ of algebraic objects
of a given type is a relation $\sigma=\tau$ in a family
of operation-symbols $g_1,g_2,\dots$ of specified arities, which holds
for every algebra $A\in\mathbf{X}$ and every assignment to each $g_i$
of a {\em derived operation} of the corresponding
arity~\cite{hyperid_WT}.

\cite[second sentence on p.65]{hyperid} \ %
If a variety ${\mathbf V}$ of groups satisfies a nontrivial
monoid identity (a monoid identity other than those satisfied
by the variety of all monoids), must ${\mathbf V}$
satisfy a hyperidentity not satisfied by the variety of all groups?
(If one merely assumes that ${\mathbf V}$ satisfies a nontrivial group
identity, the answer is no:
by \cite[Corollary~2]{hyperid} the variety of metabelian groups
satisfies the same hyperidentities as the variety of all groups.
But this is proved using the fact that
that variety satisfies no nontrivial monoid identities.)
\end{Q}

\begin{KQ}\label{Q.monoid_ids} 
Answered in \cite{SVI+AMS}.)
\end{KQ}

\begin{KQ}\label{Q.A1A2A3} 
Appears as Question 21.18 in the Kourovka Notebook.
Answered negatively.)
%
\end{KQ}

\begin{Q}\label{Q.growth} 
Question 19.110 in the Kourovka notebook asks
whether there exists a variety $V$ of groups such that
the free group on 2 generators in $V$ is finite, while the
free group on 3 generators in $V$ is infinite.

Along the same lines, one may ask,
\vspace{.3em}

(i)~~Does there exist a variety $V$ of groups such that the free group
on $2$ generators in $V$ has polynomially bounded growth, while the
free group on $3$ generators in $V$ does not?
\vspace{.3em}

The above suggests a wide class of related questions.
To state them, we need some notation regarding polynomial
bounds on growth of groups.
I set up such notation below.
If, as is likely, notation/language for this
already exists, I would appreciate a pointer to it!

Given a group $G$ with a finite generating set $X,$ let us define
the {\em growth function} $\r{gr}(G,X)$ by letting $\r{gr}(G,X)(n)$ be
the number of elements of $G$ that can be written as products
of $\leq n$ elements of $X\cup X^{-1}.$
We now define the growth degree, $\r{gr.deg}(G),$
a nonnegative real number or $+\infty,$ to be
$\r{lim\,sup}_{n\to\infty}\ \r{log}\,(\r{gr}(G,X)(n))/\r{log}\ n,$
for any choice of finite generating set $X$ for $G.$
It is not hard to show that the value so defined is independent of $X.$
Before getting to questions paralleling~(i) above, we ask,
\vspace{.3em}

(ii)~~What real numbers occur as growth degrees of
free groups $G$ on finitely many generators in group varieties?
(Must they be integers?)
\vspace{.3em}

(iii)~~Is the $\r{lim\,sup}$ in the definition of
growth degree always in fact realized as a limit; i.e., is
it equal to the corresponding $\r{lim\,inf}$?
If not always, is this true in the case of free groups
in group varieties?
\vspace{.3em}

Here, now, is a wide, though vague, generalization of question~(i):
\vspace{.3em}

(iv)~~If $V$ is a variety of groups, and for each positive
integer $i$ we write $G_i$ for the group freely generated in $V$ by
an $i$-element set, what conditions must the
sequence of values $\r{gr.deg}(G_i)$ $(i>0)$ satisfy?
\vspace{.3em}

In particular,~(i) asks whether, if $\r{gr.deg}(G_2)$ is finite,
then $\r{gr.deg}(G_3)$ must also be finite.

If the answer to~(iii) is affirmative,
one can show that the growth degree of the direct product of
two finitely generated groups is the sum of the
growth degrees of those groups.
Hence under that assumption, in the context
of~(iv), using the fact that for positive integers $i$ and $j,$ the
group $G_{i+j}$ admits a surjective homomorphism to $G_i\times G_j,$
one can deduce that
$\r{gr.deg}(G_{i+j})\geq\r{gr.deg}(G_i)+\r{gr.deg}(G_j).$
(Equality need not hold here, as can be seen
from the case $i=j=1$ for many familiar varieties of groups.)
Without knowing the answer to~(iii), one can show that
for any finitely generated group $G$ and positive integer $m,$
one has $\r{gr.deg}(G^m)= m\,\r{gr.deg}(G),$
and so deduce that for free groups in any variety, one has
$\r{gr.deg}(G_{mi})\geq m\,\r{gr.deg}(G_i).$
\end{Q}

%
%

\noindent{\bf Acknowledgements.}
I am indebted to Yves de Cornulier, Alexander Olshanskiy,
Dale Rolfsen and Dave Witte\ Morris for helpful comments on
earlier versions of this list.\vspace{.5em}

\noindent{\bf In each reference below,} the questions where
that work is cited are noted at the end, in parentheses.

\end{document}